\documentclass[journal]{IEEEtran}

\usepackage{moreverb}
\usepackage{setspace}
\usepackage{multicol}
\usepackage{float}
\usepackage{subfigure}
\usepackage{amsmath}
\usepackage{amssymb}
\usepackage{amsfonts}
\usepackage{graphicx}
\usepackage{url}
\usepackage{booktabs} 
\usepackage{color}
\usepackage{comment}
\usepackage{epstopdf}
\usepackage{indentfirst}
\usepackage{cite}

\newcommand{\diag}{{\rm diag}}

\newtheorem{Def}{Definition}
\newtheorem{Thm}{Theorem}

\newtheorem{Lem}{Lemma}
\newtheorem{Col}{Corollary}
\newtheorem{Rem}{Remark}

\newtheorem{Prm}{Problem}

\begin{document}
%
\title{Data-Based Optimal Control of
{Multi-Agent} Systems: A Reinforcement Learning Design Approach}
%
%
%

\author{Jilie~Zhang, 
        Zhanshan Wang 
        and~Hongwei Zhang 
\thanks{{This work was supported by the National Natural Science Foundation of China under grants 61603314,  61773322, 61627809,61473070 and 61433004. (Corresponding author: Hongwei Zhang)}}
\thanks{J. Zhang is with School of Information Science and Technology, and H. Zhang is with the School of Electrical Engineering, Southwest Jiaotong University, Chengdu 611756, P. R. China {(Email: jilie0226@163.com; hwzhang@swjtu.edu.cn)}; Z. Wang is with School of Information Science and Engineering, Northeastern University, Shenyang 110819, P.R. China {(Email:wangzhanshan@ise.neu.edu.cn)}.}
}

\markboth{IEEE Transactions on Cybernetics}
{Shell \MakeLowercase{J. Zhang\textit{et al.}}: Data-based optimal control of multi-Agent systems}

\maketitle

\begin{abstract}
 This paper studies optimal consensus tracking problem of heterogeneous linear multi-agent systems. {By introducing tracking} error dynamics, the optimal tracking problem is reformulated as finding a Nash-equilibrium solution of a multi-player games, which can be done by solving associated coupled Hamilton-Jacobi (HJ) equations. A data-based error estimator is designed to obtain the data-based control for the multi-agent systems.  Using the quadratic functional to approximate the every agent's value function, we can obtain the optimal cooperative control by input-output (I/O) $Q$-learning algorithm with value iteration technique in the least-square sense. The {control law} solves the optimal consensus problem for multi-agent systems with measured input-output information, and does not rely on the model of multi-agent systems. A numerical example is provided to illustrate the effectiveness of the proposed algorithm.
\end{abstract}


%
\IEEEpeerreviewmaketitle

\section{Introduction}
%
%
%
%
Originated from fields of computer science and management science, the concept of consensus has been introduced to control community to describe a collective behavior of a group of systems, called multi-agent systems. For the past two decades, consensus problem has been intensively studied in the control community, and finds its applications in vehicle formation, flocking, sensor networks, etc. For a comprehensive literature review, readers are referred to some recent survey papers \cite{cyrc,kcm,q,Oha2015Auto-survey} and monographs \cite{ren2008Book-distributed, Qu2009Book-cooperative, LewisZhang2014Book-cooperative}, and references therein.

In some scenarios, one may not only require consensus of all the agents, but also expect optimization of certain cost functions, such as the total energy. The challenge relies on the fact that each agent can only use the local information, i.e., the information of itself and its neighbors, for its controller design. This results in the so called distributed controller. Fruitful results have been achieved for distributed control of multi-agent systems notwithstanding, distributed optimal control problems are insufficiently investigated thus far.

For linear multi-agent systems, consensus tracking was guaranteed  using optimal state feedback control law derived from Riccati equations \cite{zld,hylx}; and it was shown that, using Riccati design, the control law yields an unbounded synchronization region and is robust to changes of communication graph topologies. The neural-network-based approximate dynamic programming (ADP) method has been proposed for solving the optimal synchronization of linear multi-agent systems \cite{Vamvoudakis2012,zlwf}, and for multi-player games \cite{Vamvoudakis2011,mmns2018}. The work \cite{zlwf} studies optimal output regulation problem of heterogeneous linear multi-agent systems using adaptive dynamic programming approach. An LQR-based optimal design for discrete-time linear multi-agent systems was proposed in \cite{ztll} and by using the ADP technique, the systems dynamics can be partially or completely unknown. Optimal tracking control problem of linear multi-agent systems subject to external disturbances was formulated as a multi-player zero-sum differential game, and a policy iteration algorithm was proposed to solve the corresponding coupled {Hamilton-Jacobi-Isaacs} (HJI) equations \cite{jmslv}. However, most existing works on distributed optimal control problems assume the exact knowledge of the system dynamics. By developing an ADP algorithm, \cite{jj} solves optimal control problem of a single linear system with unknown dynamics, instead of multi-agent systems.  The work \cite{zl} proposes a neural adaptive control algorithm, and solves the distributed tracking control problem of nonlinear multi-agent systems with unknown dynamics. Optimization is not considered in \cite{zl}, however.  A very recent work \cite{zzf} deals with the distributed optimal control problem of nonlinear multi-agent systems, and proposes a model free algorithm at the expense of {additional compensators}. Although \cite{ztll} allows the systems dynamics to be partially or completely unknown, all agents have {identical control matrixs}.

These above mentioned concerns motivate our research on the data-based model free optimal cooperative control for {heterogeneous} linear multi-agent systems. Inspired by the results \cite{Lewis2011,Vamvoudakis2012}, we present an on-line method to design the data-based optimal cooperative control for the consensus problem of multi-agent differential games, which brings together the game theory, data-based error estimator and input/output (I/O) $Q$-learning algorithm with value iteration. $Q$-learning is a data-based reinforcement learning technique, {which can be used to learn the solutions to optimal control problems for dynamic systems, such as \cite{Lewis2011,Tamimi2007, Liu2015}. One of the advantages of $Q$-learning is that it is able to find the available best actions without requiring the  model of the environment.} Because the output is a linear expression of the state for linear systems, we adopt {a} variant of $Q$-learning to get an ideal policy by input and output data, which is referred to as ``I/O $Q$-learning''.

For multi-agent systems, the coupled Hamilton-Jacobi (HJ) equations are set up, since each agent's control law depends on the information of itself and all its neighborhood agents. Therefore for multi-agent differential games, the optimal cooperative control relies on solving the coupled HJ equations. However it is generally difficult due to the coupling items.  Here, the I/O $Q$-learning algorithm is used to solve the problem of multi-agent differential games involving a technique known as value iteration {(VI)} \cite{gjj}. Value iteration algorithms {can be implemented in} two steps: value update and policy {improvement; and} provides an effective way of learning the solution of the coupled HJ equations. This paper mainly extends the method in \cite{Vamvoudakis2011} to deal with the optimization problem in the sense of Nash equilibrium involving the multi-agent systems with unknown dynamic. 
{Since it involves the coupling information in VI algorithm, the competitive relation is generated between agents.} To finish the competitive process and eliminate the updating data conflict, VI algorithm is implemented simultaneously for {all nodes}.

{In this paper, we first construct the tracking error dynamics and the data-based error estimator with input/output data from tracking error systems.} Similar with \cite{Lewis2011}, the estimator is used to obtain the data-based {control law} for the multi-agent systems. Then a quadratic functional is used to approximate every agent's value function at each point in the state space. Finally, we obtain the optimal cooperative control by I/O $Q$-learning algorithm with value iteration. The designed control solves the optimal consensus problem for multi-agent systems {with measured states, and} does not rely on the model of multi-agent systems.  VI algorithm is introduced to search for the optimal policies for multi-gent systems and {circumvents the difficulty} of initializing the multiple initial admissible policies at the {beginning} of policy iteration (PI) algorithms in \cite{Vamvoudakis2012} and \cite{jilie2014}.

The rest of this paper is organized as follows. Some preliminaries of graph theory and consensus problem are briefly introduced in Section \ref{Sec-preli}; the problem is formulated in Section \ref{sec3}; the optimal control law is proposed in Section \ref{Sec-control} and its numerical implementation is presented in Section \ref{sec4}; finally, a numerical example is given in Section \ref{Sec-simu} to illustrate the effectiveness of the proposed algorithm.

\vspace{-6pt}
\section{Preliminaries} \label{Sec-preli}
In this section, we briefly review the notion of consensus and some background of graph theory.

\textbf{Consensus of {Multi-Agent} Systems:}
A multi-agent system can be seen as a network which consists of a group of agents. Every agent is called a node in the network. Let $x_i \in {R}^n$  denote the state of node $v_i$. We call $\mathcal{G}_x = (\mathcal{G},x)$ a network (or algebraic graph) with the state $x \in {R}^{\bar{N} n}$, where $x=[x_1^T,\ldots,x^T_{\bar{N}}]^T$ and $\bar{N}$ is the number of nodes in the network. The state of a node may represent certain physical quantities of the agent, such as altitude, position, voltage, etc. We say all nodes reach a consensus if and only if $x_i(k)\to x_j(k)$ for all $i,j$  as $k\to \infty$. While for the consensus problem with a leader $v_0$, it requires that $x_i(k) \to x_0(k)$ for all $i$, as $k\to \infty$, where $x_0(k)$ is the state trajectory of the leader. The latter problem is also known as consensus tracking problem.

\textbf{Graph Theory:}
The topology of a communication network can be expressed by a weighted graph.
Define a weighted graph of $\bar{N}$ nodes  as $\mathcal{G} = (\mathcal{V}, \mathcal{E}, \mathcal{A})$, where $\mathcal{V}=\{{v}_1,\ldots,{v}_{\bar{N}}\}$ is the node set, $\mathcal{E}$ is the edge set with $\mathcal{E} \subseteq \mathcal{V} \times \mathcal{V}$, and $\mathcal{A}=[a_{ij}]$ is a weighted adjacency matrix. An edge of the graph $\mathcal{G}$ is denoted by $\nu_{ij}=(v_j,v_i)$, graphically represented by an arrowed link from node $j$ to node $i$, means the information can flow from node $v_j$ to node $v_i$.  In this paper, we assume that all edges are positive, i.e., $a_{ij} > 0$ if and only if
 $\nu_{ij} \in \mathcal{E}$; otherwise $a_{ij} = 0$. For notational simplicity, denote $i\in \mathcal{I}=\{1,2,\ldots,\bar{N}\}$.

\begin{Def}[Laplacian Matrix]\label{D1} The graph Laplacian matrix $L=[l_{ij}]$ is defined as $L=\mathcal{D}-\mathcal{A}$, where $\mathcal{D} = \diag \{d_i\} \in R^{\bar{N} \times \bar{N}}$ is the in-degree matrix of graph $\mathcal{G}$, with $d_i = \sum_{j=1}^{\bar{N}} a_{ij}$ being the in-degree of node $v_i$.
\end{Def}

In this paper, the simple graph is considered, i.e., no repeated edges and no self loops. Let $N_i=\{v_j\in\mathcal{V}:(v_j,v_i)\in\mathcal{E}\}$ denote the set of neighbors of node $v_i$.  A digraph is  strongly connected, if there is a directed path from node $v_i$ to node $v_j$, for all ordered pair of nodes $v_i, v_j$.  In this paper, we consider strongly connected communication graphs with fixed topology.

\section{Problem formulation}\label{sec3}
Consider the following multi-agent system with $\bar{N}$ agents
\begin{eqnarray}\label{E1}
x_i(k+1)=Ax_i(k)+B_iu_i(k),~~~i=1,2,\dots,\bar{N}
\end{eqnarray}
where $x_i(k) \in R^n$ is the state of node $v_i$, $u_i(k)\in R^{m_i}$ is the input, $A \in R^{n\times n}$  is the system matrix and $B_i \in R^{n\times m_i}$ is the control matrix. {Assume} that $(A, B_i)$ is reachable \cite{gu2012} for all $i$. In this paper, we further assume that the state $x_i$ is unmeasurable. 

The leader node dynamics is
\begin{eqnarray}\label{E2}
x_0(k+1)=Ax_0(k),
\end{eqnarray}
where $x_0(k) \in R^n$.
\vspace{0pt}

The local neighborhood tracking error $e_i(k)$ of node $v_i$ is defined as
\begin{eqnarray}\label{E3}
e_i(k)=\sum \limits_{j\in N_i}a_{ij}\left(x_i(k)-x_j(k)\right)+b_i\left(x_i(k)-x_0(k)\right),
\end{eqnarray}
where $e_i(k)=[e_{i1}(k), e_{i2}(k), \ldots, e_{in}(k)]^T \in R^n$, and $b_i$ describes the link between node $v_i$ and the leader node $v_0$. If there is a link from leader node to node $v_i$, then $b_i >0$; otherwise $b_i = 0$. We assume that $b_i>0$ for at least one $i$.
Substituting (\ref{E1}) and (\ref{E2}) {into} (\ref{E3}), we have the neighborhood tracking error dynamics
\begin{equation} \label{E4a}
e_i(k+1) = Ae_i(k) + F_iu_i(k) + \sum \limits_{j\in N_i}E_{ij} u_j(k)
\end{equation}
where $F_i=(d_i+b_i) B_i$ and $E_{ij}=-a_{ij}B_j$.  It is known that  $e_i(k)\rightarrow 0$ implies that all agents reach a consensus \cite{zl}.
We define the tracking error output as
\begin{equation}\label{E4b}
y_i(k)= C_i e_i(k)
\end{equation}
where $C_i \in R^{q_i \times n}$ is the error output matrix, and $(A,C_i)$ is observable. Note that $y_i$ is not the conventional output information of node $v_i$.


%



In this paper, the controller should not only make all the follower nodes to track the leader node, but also minimize the following local performance index for every agent
\begin{align}\label{E9}
&J_i\left(e_i(k),u_i(k),u_{(j)}(k)\right)\nonumber\\
=&\sum \limits_{l=k}^\infty r_i\left(e_i(l),u_i(l),u_{(j)}(l) \right)\nonumber \\
=&\sum \limits_{l=k}^\infty (y_i^T(l) Q_{ii} y_i(l) + u_i^T(l)
R_{ii} u_i(l) +\sum \limits_{j\in N_i} u_j^T(l)R_{ij}u_j(l)),
\end{align}
where $Q_{ii}$, $R_{ii}$, $R_{ij}$ are all positive definite matrices, and $u_{(j)}$ means the set of control laws of all neighbors of node $v_i$. The performance index  depends on the input information of node $v_i$ and its neighbors. Obviously, it is similar to multi-player games \cite{Vamvoudakis2011,Wei2011}. The policies of multi-player game have a Nash Equilibrium solution. In other words, they achieve their best interests as well as  fulfilling their common task, i.e, consensus tracking.
The Nash equilibrium  for multi-player games is defined as
\begin{Def} [Nash Equilibrium]\label{D2}
An ${N}$-tuple of control policies $\{u_1^*, u_2^*,\ldots, u_{{N}}^* \}$ is referred to as a Nash equilibrium solution for an ${N}$-player game (or multi-agent system $\mathcal{G}_x$) if
\begin{align*}
J_i^* &\triangleq  J_i(u_1^*,u_2^*,\ldots,u_i^*,\ldots,u_{{N}}^*) \leq J_i (u_1^*,u_2^*,\ldots,u_i,\ldots,u_{{N}}^*),\nonumber\\ &\qquad\qquad\qquad\qquad\qquad\qquad\qquad\qquad\quad (u_i\neq u_i^*).
\end{align*}
The ${N}$-tuple of the performance values $\{J_1^*, J_2^*, \ldots, J_{{N}}^* \}$ is known as a Nash equilibrium of the ${N}$-player game (multi-agent system $\mathcal{G}_x$).
\end{Def}

Therefore, the problem considered in this paper can be formulated as follows.
\begin{Prm} \label{P1}
Consider multi-agent systems \eqref{E1} and a leader node \eqref{E2}. Design distributed control law $u_i$, such that all follower nodes reach a consensus with the leader node, and  the local performance indices (\ref{E9}) of all nodes are minimized. In other words, we aim to find the Nash equilibrium solution subject to \eqref{E4a} and \eqref{E4b}.
\end{Prm}

\vspace{-6pt}
\section{An iterative optimal control algorithm} \label{Sec-control}

Under the given admissible control policies $u_i$ and $u_{(j)}$, define the local value function $V_i(e_i)$ for node $v_i$ as
{\begin{align}\label{E10}
&V_i\left(e_i(k)\right)\nonumber\\=&\sum \limits_{l=k}^\infty r_i\left(y_i(l),u_i(l),u_{(j)}(l)\right)\nonumber\\
=& r_i\left(y_i(k),u_i(k),u_{(j)}(k)\right) + \sum \limits_{l=k+1}^\infty r_i\left(y_i(l),u_i(l),u_{(j)}(l)\right) \nonumber\\
=&y_i(k)^T Q_{ii} y_i(k) + u_i^T(k) R_{ii} u_i(k) +\sum \limits_{j\in N_i} u_j(k)^TR_{ij}u_j(k)\nonumber \\ &+V_i\left(e_i(k+1)\right).
\end{align}}
%
The local optimal value function can be expressed as
\begin{eqnarray*}\label{E12}
V_i^*\left(e_i(k)\right) =&\min \limits_{{u}_i(l)} \sum \limits_{l=k}^\infty \Big( y_i^T(l) Q_{ii} y_i(l) + u_i^T(l) R_{ii} u_i(l)+ \nonumber\\ &\qquad\qquad\qquad\qquad \sum \limits_{j\in N_i} u_j^T(l)R_{ij}u_j(l)\Big).
\end{eqnarray*}
According to Bellman's optimality principle \cite{Zhang2013}, we have
\begin{align}\label{E13}
V_i^*\left(e_i(k)\right)=& \min \limits_{u_i(k)} \Big(y_i^T(k) Q_{ii} y_i(k) + u_i^T(k) R_{ii} u_i(k)+\nonumber\\&\ \sum \limits_{j\in N_i} u_j^T(k)R_{ij}u_j(k)\Big) + {V_i^*(e_i({k+1}))},
\end{align}
and the optimal control is given by
\begin{align}\label{E14}
u^*_i(k) =& \arg \min \limits_{u_i(k)}\big(y_i^T(k) Q_{ii} y_i(k) + u_i^T(k) R_{ii} u_i(k)+
\nonumber\\&\ \ \sum \limits_{j\in N_i} u_j^T(k)R_{ij}u_j(k)\big) + V_i^*(e_i({k+1})).
\end{align}

Then two important facts are obtained by the following Lemma \cite{Vamvoudakis2012}. 
\begin{Lem} \label{L2} Assume that the graph is strongly connected and at least one gain $b_i$
is nonzero. In addition $(A,B_i)$ is reachable for all $i$. Let $V^*_i(e_i)>0$ be a solution of the coupled HJ equation (\ref{E13}), and the optimal control $u_i^*$
be given by (\ref{E14}). Then
\begin{description}
  \item[(I)] The neighborhood tracking error system (\ref{E4a}) is asymptotically stable.
  \item[(II)] The optimal local performance $J^*_i\left(e_i(0), u^*_i, u^*_{(j)}\right)$ are equal to $V^*_i\left(e_i(0)\right)$ for all $i$, and $u_i^*, i \in{\cal I}$ are the Nash equilibrium.
\end{description}
\end{Lem}

Obviously Problem \ref{P1} is solved if the coupled HJ equations (\ref{E13}) can be solved. {However obtaining the analytical solution of \eqref{E13} is generally difficult}. In this paper, the value iteration (VI) algorithm \cite{Lewis2011} is used to solve the coupled HJ equations.

Select any set of initial control policies $u_i^0$  $(i \in \mathcal{I})$, not necessarily admissible. Then for $s = 0, 1 \ldots$, perform the following two steps until convergence:

1. Value Update:
\begin{align}\label{E16}
V_i^{s+1}\left(e_i(k)\right) &= y_i^T(k) Q_{ii} y_i(k) + u_i^{sT}(k) R_{ii} u_i^s(k) +\nonumber\\
&\sum \limits_{j\in N_i} u_j^{sT}(k)R_{ij}u_j^s(k) + V_i^{s}\left(e_i(k+1)\right).
\end{align}

2. Policy Improvement:
\begin{align}\label{E17}
u_i^{s+1}(k)=& \arg \min \limits_{u_i(k)} \big(y_i(k)^T Q_{ii} y_i(k) + u_i^{{s}}(k)^T R_{ii} u_i^{{s}}(k)+\nonumber \\ &\sum \limits_{j\in N_i}u_j^{{s}}(k)^T R_{ij}u^{{s}}_j(k)\big) + V_i^{s+1}(e_i(k+1)).
\end{align}

The coupled HJ equation (\ref{E13}) can be solved by VI algorithm iterating between (\ref{E16}) and (\ref{E17}), and the convergence of the algorithm is shown in the following Theorem.

\begin{Thm}\label{T1}  Suppose the condition $0 \leq J_i^*(e_i(k+1)) \leq \theta_i r_i(e_i, u_i, u_{(j)})$ holds uniformly for some $0 < \theta_i < \infty$. And properly pick $V_i^0 \geq 0$, $0 \leq \alpha_i \leq 1$ and $1 \leq \beta_i < \infty$, such that $0 \leq \alpha_i J_i^* \leq V_i^0 \leq \beta_iJ_i^*$ . The control sequence $\{u_i^s\}$ and value function sequence ${V_i^s}$ are iteratively updated by (\ref{E16}) and (\ref{E17}). Then the value function $V^s_i$ approaches $J_i^*$ as $s \rightarrow \infty$, according to the inequalities
\begin{align} \label{20171213Eq-1}
&\left(1+ \frac{\alpha_i-1}{(1+\theta_i^{-1})^s} \right) J^*_i(e_i(k))\nonumber\\
&\leq V_i^s(e_i(k))\leq \left(1+ \frac{\beta_i-1}{(1+\theta_i^{-1})^s} \right) J^*_i(e_i(k)).\\ \nonumber
\end{align}
\end{Thm}
\emph{Proof.}
See Appendix.
\hfill
$\Box$

\begin{Rem}
	Note that the parameters $\alpha_i$ and $\beta_i$  are only used in convergence analysis of the algorithm and have nothing to do with the controller design.
	Also the inequality $0\leq J_i^*\leq V_i^0$ holds, since $V_i^0$ is not optimal value function in general. Therefore one can easily find such $0\leq \alpha_i\leq 1$ and $1\leq \beta_i$ that $0\leq \alpha_i J_i^*\leq V_i^0\leq \beta_i J_i^*$.
\end{Rem}

According to the results of Theorem \ref{T1}, we have the following corollary.

\begin{Col} Define the control sequence $\{u_i^s\}$ as in (\ref{E17}) and the value function $\{V_i^s\}$ as in (\ref{E16}). When $\{V_i^s\}\rightarrow J^*_i$  as $s \rightarrow \infty$, the control sequence ${u_i^s}$ converges to the optimal control law $u_i^*$  as $s \rightarrow \infty$, that is, $\lim \limits_{s \rightarrow \infty} u_i^s(k) = u_i^*$.
\end{Col}

The coupled HJ equation (\ref{E13}) is constructed by the value function $V_i(e_i(k))$ (\ref{E10}) subject to the tracking error dynamics (\ref{E4a}) and (\ref{E4b}). Once the optimal solution (the value function $V_i^*(e_i(k))$) is obtained, the optimal admissible control (\ref{E14}) based on $V_i^*(e_i(k))$ makes the tracking error $e_i(k)$ tend to zero as $k \rightarrow\infty$. In other words, all follower nodes track the leader node as $k \rightarrow\infty$. From above analysis, it is clearly seen that both tracking and optimality of the multi-agent systems are obtained simultaneously.

\vspace{-6pt}
\section{Data-based implementation of the optimal control algorithm}\label{sec4}
In Subsection \ref{sec3-sub1}, we construct an error estimator using the measured input/output data from tracking error dynamics \eqref{E4a} and \eqref{E4b}. Then in Subsection  \ref{sec3-sub2} the nearly data-based cooperative control is designed by the I/O $Q$-learning method \cite{Tamimi2007}. Subsection \ref{sec3-sub3} presents how to solve the optimal inner kernel matrix $\bar{P}_i$ by VI algorithm.

\vspace{-2pt}
\subsection{Data-based Error Estimator}\label{sec3-sub1}
Hereinafter we derive a tracking error estimator.
Given the current time $k$, the tracking error dynamics can be written as the expanded error equation on the time horizon $[k-N, k]$,
\begin{subequations}\label{E18}
\begin{align}
\label{E18a}
e_i(k)=& {A}^N {e}_i(k-N) + \mathcal{B}_{N_i} \bar{u}_{i[k-1,k-N]}\nonumber \\
& +{\sum \limits_{j\in{N}_i} \mathcal{B}_{N_{ij}}\bar{u}_{j[k-1,k-N]}}\\
\label{E18b}
\bar{y}_{i[k-1,k-N]} =& \mathcal{C}_{N_i} {e}_i{(k-N)} + \mathcal{D}_{N_i} \bar{u}_{i[k-1,k-N]}\nonumber \\
 &+ {\sum \limits_{j\in{N}_i}\mathcal{D}_{N_{ij}} \bar{u}_{j[k-1,k-N]}}
\end{align}
\end{subequations}
where
\begin{align*}
&\bar{u}_{i[k-1,k-N]}= [u_i({k-1}), u_i({k-2}), \cdots, u_i({k-N})] \in {R}^{m_i N},\\
& \bar{u}_{j[k-1,k-N]}=[u_j({k-1}), u_j({k-2}),  \cdots, u_j({k-N})] \in {R}^{m_j N},\\
&\bar{y}_{i[k-1,k-N]}= [y_i({k-1}),   y_i({k-2}),  \cdots,  y_i({k-N})] \in {R}^{q_i N},\\
&\mathcal{B}_{N_i} = [F_i \ AF_i\ \cdots\ A^{N-2}F_i\  A^{N-1}F_i]\in R^{n\times m_iN},\\
&\mathcal{B}_{N_{ij}} = [E_{ij}\ AE_{ij}\ \cdots\  A^{N-2}E_{ij}\  A^{N-1}E_{ij}]\in R^{n\times m_jN}, 
\end{align*}
\begin{align*}
\mathcal{C}_{N_i}  =& \left[
                     \begin{array}{c}
                        C_iA^{N-1} \\
                        C_iA^{N-2} \\
                        \vdots \\
                        C_i \\
                     \end{array}
                   \right]\in {R}^{q_i N\times n},\\
\mathcal{D}_{N_i} =& \left[
                     \begin{array}{ccccc}
                       0 &~ C_iF_i &~ \cdots &~ C_iA^{N-3}F_i &~ C_iA^{N-2}F_i \\
                       0 &~ 0 &~ \cdots &~ C_iA^{N-4}F_i &~ C_iA^{N-3}F_i \\
                       0 &~ 0 &~ \ddots &~ \vdots &~ \vdots \\
                       0 &~ 0 &~ \cdots &~ 0 &~ C_iF_i \\
                       0 &~ 0 &~ \cdots &~ 0 &~ 0\\
                     \end{array}
                   \right]\\
&\qquad \qquad \qquad \qquad \qquad \qquad \qquad \quad\quad \in {R}^{q_i N\times m_iN},\\
\mathcal{D}_{N_{ij}}=& \left[
          \begin{array}{ccccc}
              0 &~C_iE_{ij} &~\cdots &~C_iA^{N-3}E_{ij} &~C_iA^{N-2}E_{ij} \\
              0 &~0  &~\cdots &~C_iA^{N-4}E_{ij} &~C_iA^{N-3}E_{ij} \\
         \vdots &~\vdots &~\ddots &~\vdots&~\vdots\\
              0 &~0 &~\cdots &~0 &~C_iE_{ij}\\
              0 &~0 &~\cdots &~0 &~0 \\
          \end{array}
        \right]\\ &\qquad \qquad \qquad \qquad \qquad \quad \quad \quad \in {R}^{q_iN\times m_jN}, j \in N_i.
\end{align*}

In (\ref{E18}), vectors $\bar{u}_{i[k-1,k-N]}$ and $\bar{u}_{j[k-1,k-N]}$ are the input sequences over the time interval $[k-N,k-1]$, while the vector $\bar{y}_{i[k-1,k-N]}$ is output sequence over the time interval $[k-N,k-1]$. They denote the available measured data. $\mathcal{B}_{N_i}$ and $\mathcal{B}_{N_{ij}}$ are controllability matrices of ($A, F_i$) and ($A, E_{ij}$). $\mathcal{D}_{N_i}$ and $\mathcal{D}_{N_{ij}}$ represent the relation between $\bar{u}_{i[k-1,k-N]}$, $\bar{u}_{j[k-1,k-N]}$ and
$\bar{y}_{i[k-1,k-N]}$. $\mathcal{C}_{N_i}$ is the observability matrix of $(A,C_i)$.

Note that since $(A,C_i)$ is observable, there exists an observability index $K_i$ such that $rank(\mathcal{C}_{N_i}) < n$ for $N <K_i$ and that $rank(\mathcal{C}_{N_i})=n$ for $N \geq K_i$. Note that $K_i$ satisfies $K_iq_i \geq n$. Let $N \geq K_i$. Therefore $\mathcal{C}_{N_i}$ has full column rank of $n$, namely there exists a left inverse for $\mathcal{C}_{N_i}$. If $N$ is larger than $n$, then there exists a left inverse for $\mathcal{C}_{N_i}$ due to the observability of $(A,C_{i})$, and the observability index satisfying $K_{i}\leq n$. In general the rough information should be known, even if the model of system is unknown. If only the dimension of system $n$ is known or $N$ is selected as the sufficiently large constant, the estimator can be established. Furthermore we can obtain Theorem \ref{T2}.
\begin{Thm}\label{T2}
If $C_{N_i}$ is of full column rank, then there exists an error estimator based on the input/output data sequences $\bar{u}_{i[k-1,k-N]}$, $\bar{u}_{j[k-1,k-N]}$ and $\bar{y}_{i[k-1,k-N]}$ for the tracking error system (\ref{E4a}) and \eqref{E4b}, as follows
\begin{align}\label{E19}
e_i(k)=& {T}_{u_i} \bar{u}_{i[k-1,k-N]}+\sum \limits_{j\in{N}_i} {T}_{u_{ij}} \bar{u}_{j[k-1,k-N]}\nonumber \\
 &\qquad\qquad\qquad\qquad+{T}_{y_i} \bar{y}_{i[k-1,k-N]},
\end{align}
where ${T}_{y_i}$ = ${A}^N\mathcal{C}^+_{N_i}$, ${T}_{u_i}$ = $\mathcal{B}_{N_i}- {A}^N\mathcal{C}_{N_i}^+\mathcal{D}_{N_i}$ and ${T}_{u_{ij}}$ = $\mathcal{B}_{N_{ij}}-{A}^N\mathcal{C}^+_{N_i}\mathcal{D}_{N_{ij}}$ with $\mathcal{C}_{N_i}^+=(\mathcal{C}_{N_i}^T\mathcal{C}_{N_i})^{-1}\mathcal{C}_{N_i}^T$ being the left pseudo inverse of
$\mathcal{C}_{N_i}$.
\end{Thm}
\emph{Proof.}
It is similar to the proof of Lemma 1 in \cite{Lewis2011}. Thus details are omitted.  \hfill $\Box$

\begin{Rem}\label{R1}
From (\ref{E19}) we can see that the error variable is completely eliminated from the right side of (\ref{E4a}). That is, $e_i(k)$ can be reconstructed by the available measured input/output data sequences $\bar{u}_{i[k-1,k-N]}$, $\bar{u}_{j[k-1,k-N]}$ and $\bar{y}_{i[k-1,k-N]}$. Therefore (\ref{E19}) is referred to as a data-based error estimator for the tracking error system (\ref{E4a}).
\end{Rem}
\vspace{-2pt}
\subsection{The Implementation of Data-based Control}\label{sec3-sub2}
In this subsection, according to Theorem \ref{T2}, we obtain the nearly data-based optimal cooperative control $u_i(k)$ by I/O $Q$-learning technology,
using measured input/output data from the tracking error dynamic systems (\ref{E4a}) and \eqref{E4b}.

First, let {$\bar{w}_{i[k-1,k-N]} = [\bar{U}_{{i[k-1,k-N]}}, \quad \bar{y}_{i[k-1,k-N]}], $} where $\bar{U}_{{i[k-1,k-N]}} \triangleq [\bar{u}_{{i[k-1,k-N]}}, \bar{U}_{{j[k-1,k-N]}}]$ with $\bar{U}_{{j[k-1,k-N]}}$  being the column vector of $ \{\bar{u}_{{j[k-1,k-N]}} | j \in {N}_i\}$. Then (\ref{E19}) can be written as
\begin{equation}\label{E22}
e_i(k)=\left[
                       \begin{array}{ccc}
                         {T}_{u_i} &{ \bar{T}_{u_{ij}}} & {T}_{y_i} \\
                       \end{array}
                     \right]
\bar{w}_{i[k-1,k-N]},
\end{equation}
where ${\bar{T}_{u_{ij}}}\triangleq [{T}_{u_{i1}}\quad \ldots \quad {T}_{u_{ij}}\quad \ldots \quad {T}_{u_{iN}} ], j \in \mathcal{I}, j\neq i$.

For the tracking error systems (\ref{E4a}), although the explicit form of the actual cost functional $V_i(e_i(k))$ is unclear, it can be approximated by a quadratic functional at each point in the state space (similarly, see \cite{Murry2002}). So the approximated cost functional associated to any policy ${u}_i(k)$ and ${u}_{(j)}(k)$ (not necessarily optimal) is
\begin{equation}\label{E23}
V_i(e_i(k)) = e_i^T(k) {P_i} e_i(k),
\end{equation}
where $P_i$ is a symmetric positive definite matrix.

We can obtain the approximation cost functional by (\ref{E22}) and (\ref{E23}) as
\begin{align}\label{E24}
&V_i(e_i(k))\nonumber\\
=&e_i^T(k) {P_i} e_i(k) \nonumber\\
=& \bar{w}_{i[k-1,k-N]}^T   \left[
                                                 \begin{array}{c}
                                                   {T}_{u_i}^T \\
                                                   {\bar{T}_{u_{ij}}^T} \\
                                                   {T}_{y_i}^T \\
                                                 \end{array}
                                               \right]
 P_i  \left[
                       \begin{array}{ccc}
                         {T}_{u_i} &{ \bar{T}_{u_{ij}}} & {T}_{y_i} \\
                       \end{array}
                     \right]
\bar{w}_{i[k-1,k-N]}\nonumber\\
=& \bar{w}_{i[k-1,k-N]} ^T \bar{P}_i
 \bar{w}_{i[k-1,k-N]},
\end{align}
where $\bar{P}_i=\left[
  \begin{array}{ccc}
    {T}_{u_i}^T P_i {T}_{u_i} & {T}_{u_i}^T P_i {\bar{T}_{u_{ij}}} & {T}_{u_i}^T P_i {T}_{y_i} \\
    {\bar{T}_{u_{ij}}^T} P_i {T}_{u_i} & {\bar{T}_{u_{ij}}^T} P_i {\bar{T}_{u_{ij}}} & {\bar{T}_{u_{ij}}^T} P_i {T}_{y_i}  \\
    {T}_{y_i}^T P_i {T}_{u_i} & {T}_{y_i}^T P_i {\bar{T}_{u_{ij}}} & {T}_{y_i}^T P_i {T}_{y_i}  \\
  \end{array}
\right]$.
%
%
Substituting (\ref{E24}) into (\ref{E10}), we have
\begin{align}\label{E26}
&V_i(e_i(k)) \nonumber \\
=& \bar{w}_{i[k-1,k-N]} ^T \bar{P}_i \bar{w}_{i[k-1,k-N]} \nonumber \\
=& y_i^T(k) Q_{ii} y_i(k) + u_i^T(k) R_{ii} u_i(k) +\sum \limits_{j\in {N}_i} u_j^T(k)R_{ij}u_j(k)\nonumber \\
& + \bar{w}_{i[k,k-N+1]} ^T \bar{P}_i \bar{w}_{i[k,k-N+1]}.
\end{align}
Partition $\bar{w}_{i[k,k-N+1]} ^T\bar{P}_i \bar{w}_{i[k,k-N+1]}$ as
\begin{align}\label{E27}
&\bar{w}_{i[k,k-N+1]}^T\bar{P}_i \bar{w}_{i[k,k-N+1]}  \nonumber \\
=& \left[
             \begin{array}{c}
               u_i(k) \\
               \bar{u}_{i[k-1,k-N+1]} \\
               {\bar{U}_{j[k,k-N+1]}} \\
               \bar{y}_{i[k,k-N+1]} \\
             \end{array}
           \right]^T
           \bar{P}_i
            \left[
             \begin{array}{c}
                u_i(k) \\
               \bar{u}_{i[k-1,k-N+1]} \\
               {\bar{U}_{j[k,k-N+1]}} \\
               \bar{y}_{i[k,k-N+1]} \\
             \end{array}
           \right],
\end{align}
where

\def\hsymb#1{\mbox{\strut\rlap{\smash{\Huge$#1$}}\quad}}
$$\bar{P}_i= \left[
             \begin{array}{cccc}
p_{u_i u_i}            &~ p_{u_i \bar{u}_i}                 &~ \bar{p}_{u_i {\bar{U}_j}}      &~ p_{u_i \bar{y}_i}    \\
p_{u_i \bar{u}_i}^T    &~ p_{\bar{u}_i \bar{u}_i}           &~ \bar{p}_{\bar{u}_i \bar{U}_j}  &~ p_{\bar{u}_i \bar{y}_i}\\
\bar{p}_{u_i \bar{U}_{j}}^T   &~ \bar{p}_{\bar{u}_i \bar{U}_{j}}^T &~ \bar{p}_{\bar{U}_j \bar{U}_{j}} &~ \bar{p}_{\bar{U}_{j} \bar{y}_i}\\
p_{u_i \bar{y}_i }^T          &~ p_{\bar{u}_i \bar{y}_i}^T         &~ \bar{p}_{\bar{U}_j \bar{y}_i}^T  &~ p_{\bar{y}_i \bar{y}_i}\\
             \end{array}
           \right],$$
with $\bar{p}_{u_i {\bar{U}_j}} \triangleq [{p}_{u_i {\bar{u}_j}}]$, $\bar{p}_{\bar{u}_i {\bar{U}_j}} \triangleq [{p}_{\bar{u}_i {\bar{u}_j}}]$, $\bar{p}_{\bar{U}_j {\bar{U}_j}} \triangleq [{p}_{\bar{u}_j {\bar{u}_j}}]$ and $\bar{p}_{\bar{U}_j {\bar{y}_i}} \triangleq [{p}_{\bar{u}_j {\bar{y}_i}}]$ being row vectors of dimension $|N_i|$, where $j\in N_i$ and $|N_i|$ is the cardinality of $N_i$.

According to the necessary condition of optimality principle {$$\frac{\partial{V_i(e_i(k))}}{\partial{u}_i(k)} = 0,$$}
we can obtain the optimal control law
\begin{align}\label{E28}
u_i(k) =& -(R_{ii} + p_{u_i u_i})^{-1} (p_{u_i \bar{u}_i} \bar{u}_{i[k-1,k-N+1]}  \nonumber \\
&+ \bar{p}_{u_i \bar{U}_j} \bar{U}_{j[k,k-N+1]} + p_{u_i \bar{y}_i} \bar{y}_{i[k,k-N+1]})
\end{align}
where
$p_{u_i u_i} $=$ \mathcal{H} \bar{P}_i \mathcal{H}^T \in {R}^{m_i\times m_i}$, $p_{u_i \bar{u}_i}$ = $\mathcal{S} \bar{P}_i$ $\in {R}^{m_i\times m_i(N-1)}$, $ \bar{p}_{u_i \bar{U}_j} = \mathcal{W}\bar{P}_i$ $\in {R}^{m_i\times \sum \limits_{j \in N_i} m_jN}$ and $p_{u_i \bar{y}_i}$ = $\mathcal{L} \bar{P}_i \in {R}^{m_i\times q_i N}$, with $\mathcal{H}$ = $[I_{m_i} ~ 0_{m_i\times ((\sum \limits_{j \in N_i}m_j + q_i )N-m_i)}]$; $\mathcal{S}$ being an operator that yields the matrix $\left(m_i \times m_i(N-1)\right)$, which is $1_{st}$ to $m_{i_{th}}$ rows and $\left(m_i+1\right)_{th}$ to $\left({m_iN}\right)_{th}$ columns of matrix $\bar{P}_i$; $\mathcal{W}$ being an operator that yields the matrix $\left(m_i \times \sum \limits_{j \in N_i} m_jN\right)$, which is $1_{st}$ to $m_{i_{th}}$ rows and $\left(m_iN+1\right)_{th}$ to $\left(\sum \limits_{j \in \{N_i, i\}} m_jN\right)_{th}$ columns of matrix $\bar{P}_i$ and $\mathcal{L}$ being an operator that yields the matrix ($m_i \times q_i N$), which is $1_{st}$ to $m_{i_{th}}$ rows and $\left(\sum \limits_{j \in \{N_i, i\}} m_jN+1\right)_{th}$ to $\left(\sum \limits_{j\in \{N_i, i\}}(m_j+q_i)N\right)_{th}$ columns of matrix $\bar{P}_i$.
Obviously, the control law is  distributed in the sense that it only involves the neighborhood information.

The cooperative control ${u}_i(k)$ depends on $ p_{u_i u_i}, p_{u_i \bar{u}_i}, \bar{p}_{u_i \bar{U}_j}, p_{u_i \bar{y}_i}$ and the measured input/output data $\bar{u}_{i[k-1,k-N+1]}, \bar{u}_{j[k-1,k-N+1]}, \bar{y}_{i[k,k-N+1]}$. Therefore the optimal control ${u}_i(k)$ can be obtained, if only we know the ideal inner kernel $\bar{P}_i$. In the following section, we shall show how to optimize it.
\vspace{-2pt}
\subsection{Obtaining the Local Optimal Inner Matrix ${\bar{P}}_i$ by Value Iteration Algorithm}\label{sec3-sub3}
As shown in \cite{Tamimi2007}, since the control input appears in quadratic form (\ref{E10}), the minimization of (\ref{E10}) can be carried out in terms of the learned inner kernel matrix $\bar{P}_i$ without any knowledge of the system (\ref{E4a}).  This is a data-based  control approach, since the control (\ref{E28}) contains present and past observations of the input and output information.

Inspired by  \cite{Lewis2011, Tamimi2007}, we use the following VI algorithm to calculate the ideal $\bar{P}_i$ by applying only the measured input and output data, without knowledge of the state-space model of tacking error systems.

\textbf{VI algorithm for multi-agent systems:}

Select a set of control policies $u^0_i = \mu^0_i ({i=1,\ldots, \bar{N}})$. For $s = 0, 1, \ldots,$ perform the following iteration until convergence: \\
Step 1. Value Update: \\
By (\ref{E26}), solve $\bar{P}_i^{s+1}$ from
\begin{align}\label{E32}
vec(\bar{P}^{s+1}_i)=&{(\bar{w}_{i[k-1,k-N]} \otimes \bar{w}_{i[k-1,k-N]} ) ^{-1}}(y_i^T(k) Q_{ii} y_i(k)  \nonumber \\&+ u_i^{sT}(k) R_{ii} u_i^s(k)+\sum \limits_{j\in N_i} u_j^{sT}(k)R_{ij}u_j^s(k)  \nonumber \\ &+{\bar{w}_{i[k,k-N+1]}}^T\bar{P}^{s}_i \bar{w}_{i[k,k-N+1]}),{i=1,\ldots, \bar{N}}.
\end{align}
%
Step 2. Policy Improvement: \\
According to (\ref{E28}), define the updated policy by
\begin{align}\label{E30}
u_i^{s+1}(k) =& -(R_{ii} + p_{u_i u_i}^{s+1})^{-1} (p_{u_i\bar{u}_i}^{s+1} \bar{u}{^s_{i[k-1,k-N+1]}}  \nonumber \\
&+ \bar{p}_{u_i \bar{U}_j}^{s+1} \bar{U}{^s_{j[k,k-N+1]}}+ p_{u_i \bar{y}_i}^{s+1} \bar{y}{^s_{i[k,k-N+1]}}), \nonumber \\&\qquad\qquad\qquad\qquad\qquad\qquad\qquad{i = 1,\ldots, \bar{N}}.
\end{align}
Go to step 1. Stop the iteration until $u_i$ converge to $u^*_i$ (i.e., {$\|\bar{P}_i^{s+1}-\bar{P}^s_i\|\leq\varepsilon$} for an ideal constant $\varepsilon$ ) for all $i$.


\begin{Rem}
  Here we extend the method in \cite{Lewis2011} to deal with multi-agent systems. What we consider is the optimization problem in the sense of Nash equilibrium of multi-agent systems, rather than the optimization problem of the single plant as in \cite{Lewis2011}. Furthermore the implementation of our algorithm is also different from that in \cite{Lewis2011}, since it involves the coupling information in VI algorithms, see ((\ref{E32}) and (\ref{E30})) and the equation ((\ref{E3}) and (\ref{E4a})). Each agent is not considered independently. They are coupled in the tracking error $e_i(k)$ involving the information of neighbours (such as (\ref{E3}) and (\ref{E4a})). In fact, the process that VI algorithm is implemented simultaneously for all nodes can be seen as the game process between agents.
\end{Rem}

\begin{figure}[!htb]
\centering
\includegraphics[width=3.0in]{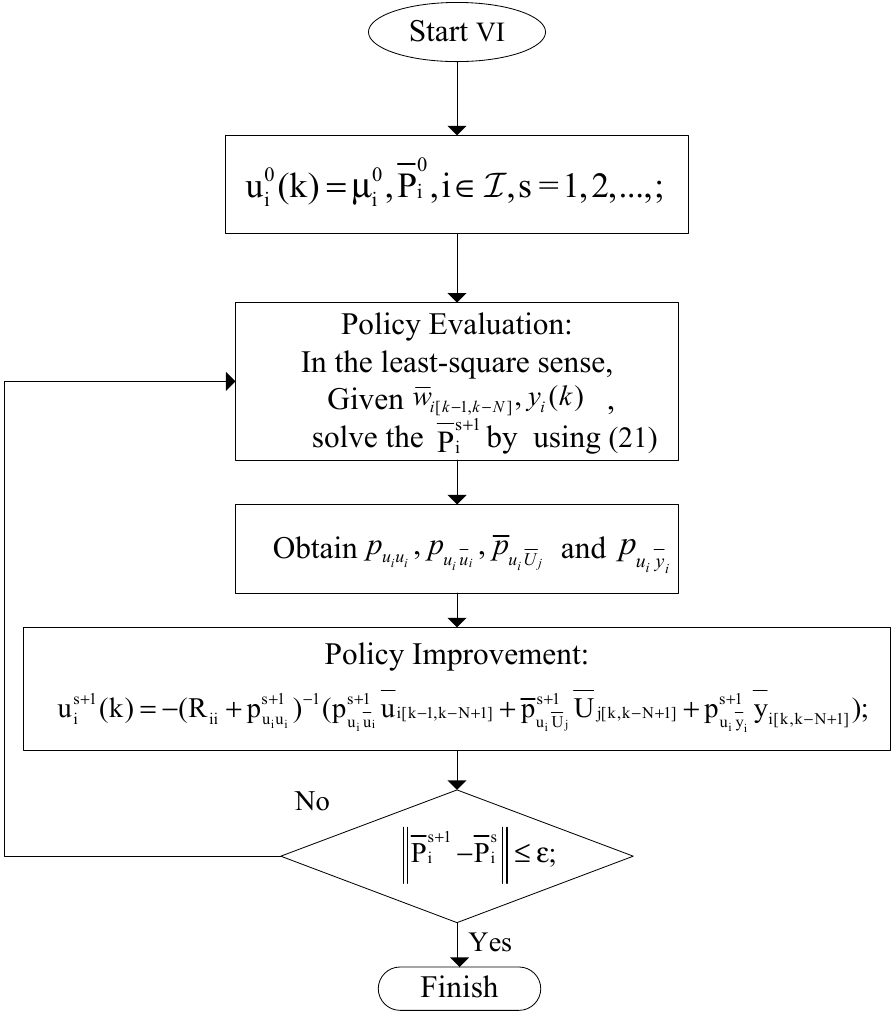}
\caption{Flowchart of value iteration algorithm} \label{F1}
\end{figure}

Solving equation (\ref{E32}) only requires the input/output data from tracking error systems. The solution of equations (i.e., the kernel matrix $\bar{P}_i \in {R}^{(\sum \limits_{j\in\{N_i, i\}}m_j + q_i)N \times (\sum \limits_{j\in \{N_i, i\}}m_j + q_i )N}$) is symmetric and has $(\sum \limits_{j\in \{N_i, i\}}m_j + q_i )N ((\sum \limits_{j\in \{N_i, i\}}m_j + q_i )N+1)/2$ independent terms. Therefore, it is necessary to sample data for at least $(\sum \limits_{j\in \{N_i, i\}}m_j + q_i )N ((\sum \limits_{j\in \{N_i, i\}}m_j + q_i )N+1)/2$ time steps for a solution using RLS. The flowchart of the VI algorithm is shown in Fig. \ref{F1}.

\begin{Rem}\label{R6} To update $\bar{P}_i$ by VI algorithm, it is necessary to insert a persistent excitation (PE) into the  control input of each agent. For details, please refer to Remark 8 in \cite{jilie2014}.
\end{Rem}

Next, we will give the detailed design procedure.

\textbf{Step 1:} Set up the local neighborhood tracking error dynamics (\ref{E4a}) by the definition of the local neighborhood tracking error (\ref{E3}).

\textbf{Step 2:} Construct the data-based error estimators (\ref{E19}) by the method \cite{Lewis2011}.

\textbf{Step 3:} Use the quadratic forms to approximate the
solution (value functions), such as (\ref{E23}) or (\ref{E24}). Substituting (\ref{E23}) into (\ref{E10}), then the optimal control law (\ref{E28}) can be obtained by solving the necessary condition of optimal principle.

\textbf{Step 4:} The optimal control can be obtained by VI algorithm for the consensus tracking problem of multi-agent systems. See Fig. \ref{F1}.

\vspace{-6pt}
\section{Numerical Example} \label{Sec-simu}

\begin{figure}[!htb]
\centering
\includegraphics[width=2 in]{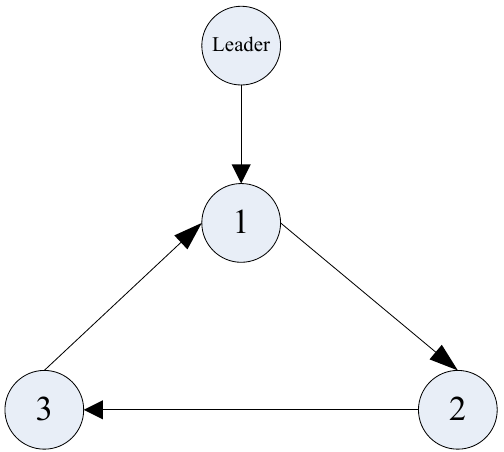}
\caption{The graph topologyx}\label{F2}
\end{figure}

Here we consider the three-node digraph with leader node
connected to node 1, shown in Fig. \ref{F2}. The edge weights are chosen as one.

The node dynamics are {
\begin{align*}
{x}_1(k+1)=& A x_1(k) + B_1u_1(k),\\
{x}_2(k+1)=& A x_2(k) + B_2u_2(k),\\
{x}_3(k+1)=& A x_3(k) + B_3u_3(k)
\end{align*}}
and the leader node is {${x}_0(k+1)=Ax_0(k)$}, where $A=\left[
                                                                      \begin{array}{cc}
                                                                        0 & 1 \\
                                                                        -1 & 0 \\
                                                                      \end{array}
                                                                    \right]$,
$B_1=\left[
         \begin{array}{c}
           2 \\
           1 \\
         \end{array}
     \right],
B_2=\left[
         \begin{array}{c}
           2 \\
           3 \\
         \end{array}
\right],
B_3=\left[
         \begin{array}{c}
           2 \\
           2 \\
         \end{array}
\right]$.
For $i=1,2,3$, let $Q_{ii}=I$; $R_{ii}=2$; $R_{ij}=0.1 (i\neq j)$ (Note that
$R_{ij}=0$, if ${v}_j$ is not the neighbor of ${v}_i$) and
$a_i=0.1$.

Our objective is to design the local data-based control law,
make $x_i$ reach a consensus with the leader and at the same time minimize the cost
functional (\ref{E9}). Letting $N=2$ and $C_i=I_2$, we can obtain the optimal
cooperative control by (\ref{E28})
\begin{align}\label{E33}
 u_i^* =& -(R_{ii} + p^*_{u_i u_i})^{-1} (p^*_{u_i \bar{u}_i} \bar{u}_{i[k-1,k-1]}   \notag \\
        & + \bar{p}^*_{u_i \bar{U}_j} \bar{U}_{j[k,k-1]} + p^*_{u_i \bar{y}_i} \bar{y}_{i[k,k-1]}).
\end{align}
\begin{figure}[!htb]
\centering
\begin{center}
\subfigure[]{
\resizebox*{6.8cm}{!}{\includegraphics{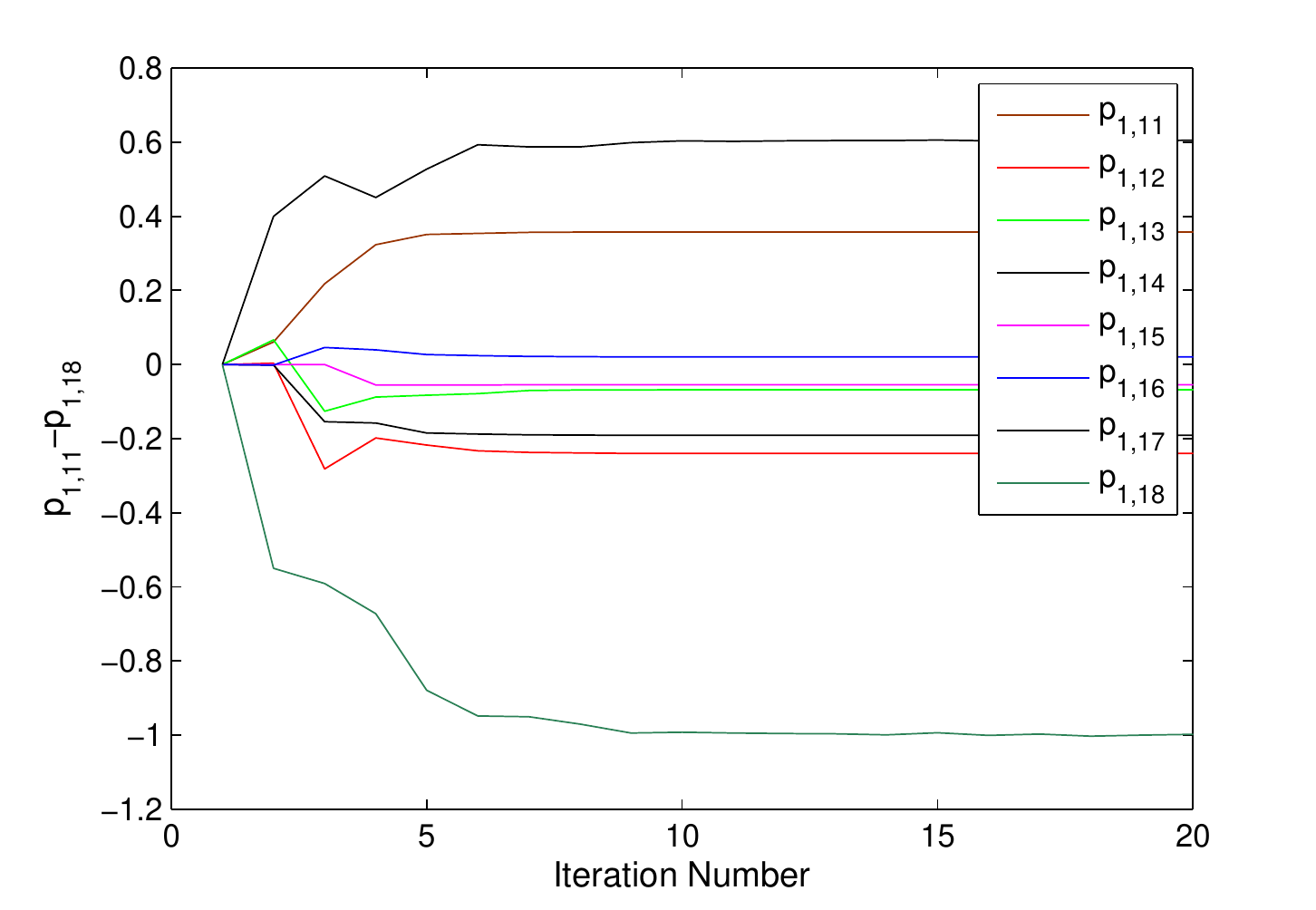}}}\label{F3}
\subfigure[]{
\resizebox*{6.8cm}{!}{\includegraphics{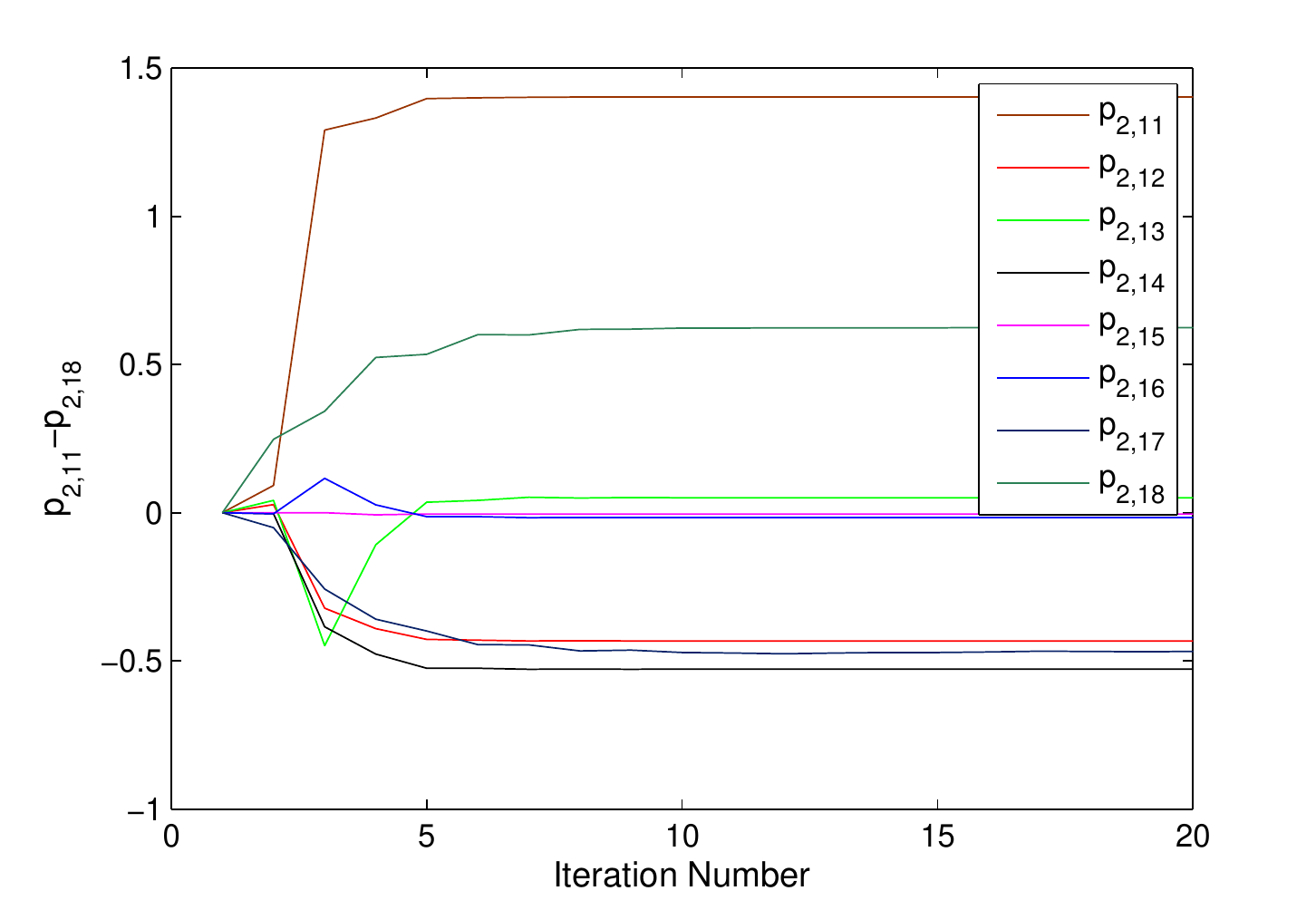}}}
\subfigure[]{
\resizebox*{6.8cm}{!}{\includegraphics{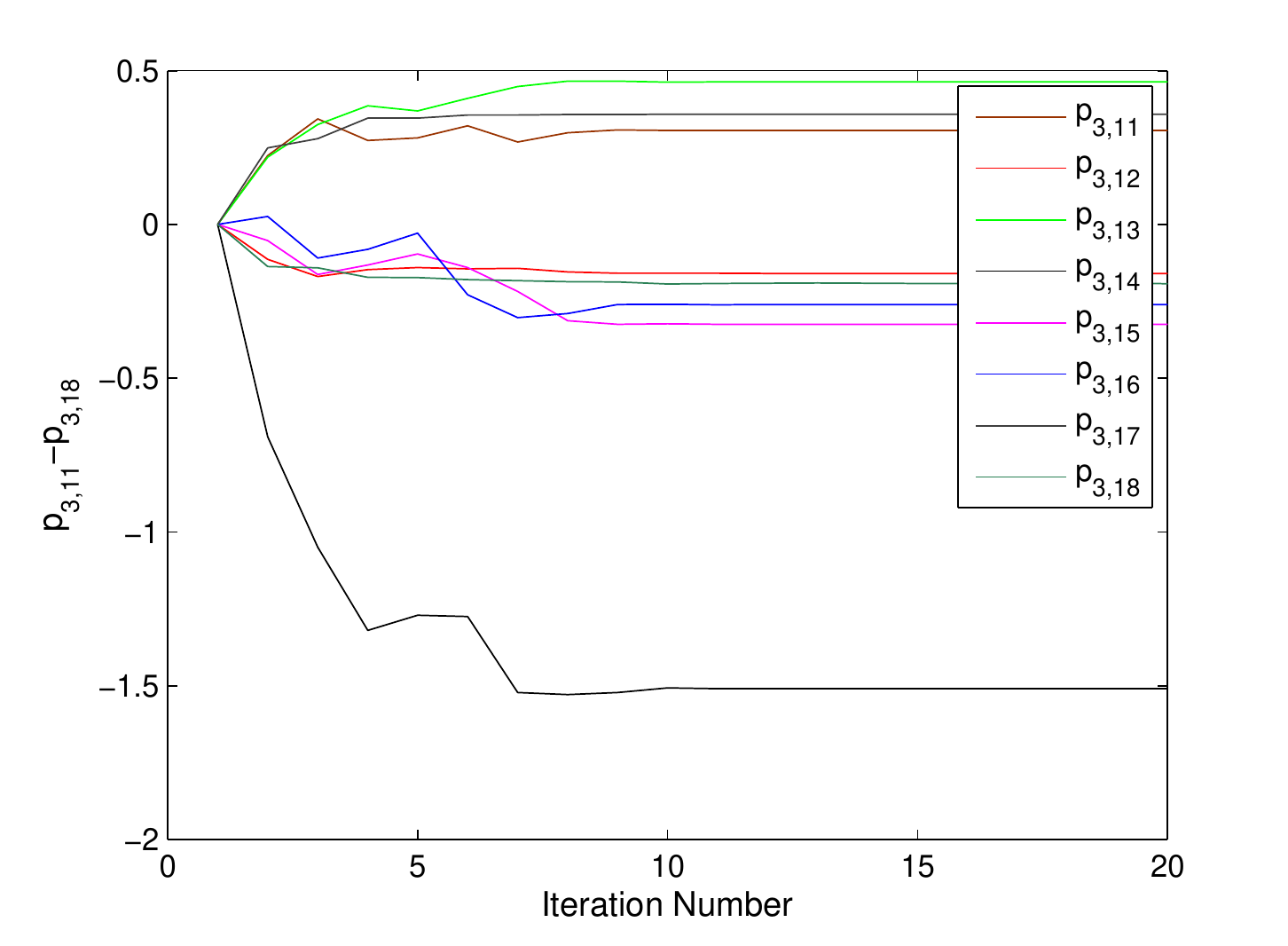}}}%
\caption{The updating profiles of weights for every agents.}%
\end{center}
\end{figure}

\begin{figure}[!htb]
\centering
\includegraphics[width=3in]{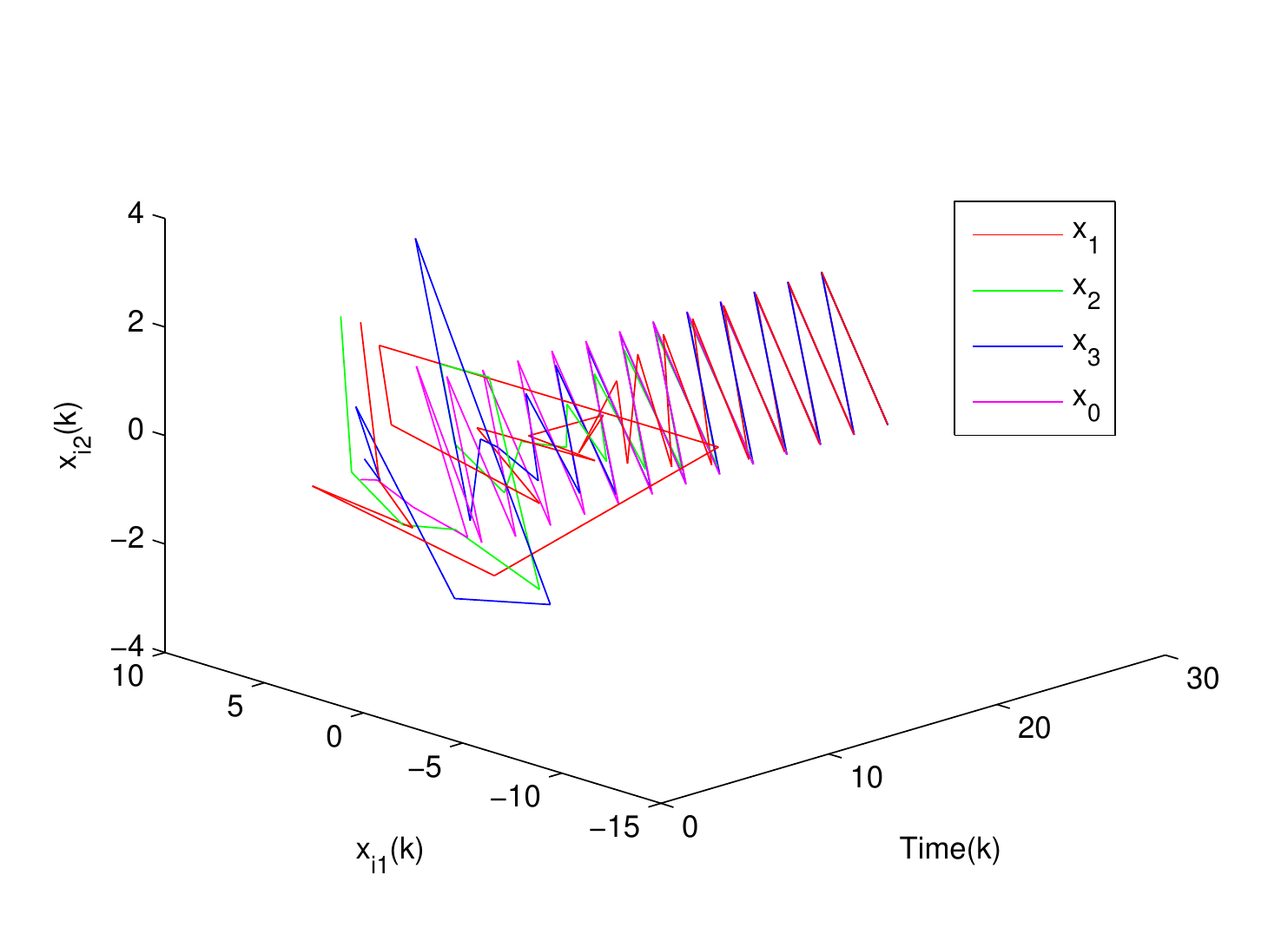}
\caption{The evolution of the system state}\label{F4}
\end{figure}

The evolution of weight $p_{i,{11}}$ to $p_{i,{18}}$ $(i=1,2,3)$ with the random initial value are shown in Fig. 3.
Obviously, after 13 iteration steps, $p_{i,{11}}$ to $p_{i,{18}}$ converge to the
ideal values.

Figure \ref{F4} depicts the evolution of the system state under the
optimal cooperative controls (\ref{E33}) with the obtained ideal
value $p_{i,{11}}$ to $p_{i,{18}}$. After 25 steps, the state of each node reaches a
consensus with leader node.

 Compared with other conventional algorithms, such as \cite{Vamvoudakis2012} and \cite{jilie2014}, this proposed algorithm features independence of the agent dynamics. In addition, this algorithm uses the value iteration to search for the optimal control law, rather than policy iteration. By doing this, we do not require the initial policy to be an admissible policy. Thus the implementation of the algorithm is significantly eased.

\vspace{-6pt}
\section{Conclusion}
A data-based optimal cooperative control law has been proposed for multi-agent systems, which brings together game theory, data-based tracking error estimator design and I/O $Q$-learning algorithm with value iteration. The key is to approximate the solution of the coupled Hamilton-Jacobi (HJ) equation, using the input/output data under the framework of the VI algorithm. Then the approximated solution is utilized to obtain the data-based optimal cooperative control law. Therefore, the proposed method designs a nearly optimal cooperative control without knowing the model of multi-agent systems. An example illustrates the effectiveness of the proposed algorithm. 

\onecolumn
\section*{ APPENDIX}
\vspace{-2pt}
\section*{Proof of Theorem \ref{T1}}
The first part of the inequality \eqref{20171213Eq-1}
\begin{eqnarray*}\label{E34}
\left(1+ \frac{\alpha_i-1}{(1+\theta_i^{-1})^s} \right) J^*_i(e_i(k))\leq V_i^s(e_i(k))
\end{eqnarray*}
will be proven by mathematical induction. When $s = 1$, since
\begin{eqnarray}\label{E35}
\frac{\alpha_i-1}{(1+\theta_i)}{\left\{\left(\theta_i r_i(e_i(k), u_i, u_{(j)})\right)-J^*_i(e_i(k+1))\right\}} \leq 0,\quad 0\leq\alpha_i\leq 1
\end{eqnarray}
and $\alpha_i J_i^* \leq ≤ V^0_i, \forall e_k(k)$ , we have
\begin{align}\label{E36}
V_i^1(e_i(k))=&\min\limits_{u_i(k)} \left\{r_i(e_i(k), u_i, u_{(j)}) + V^0_i(e_i(k+1))\right\}\nonumber\\
\geq& \min\limits_{u_i(k)} \left\{r_i(e_i(k), u_i, u_{(j)}) + \alpha_iJ^*_i(e_i(k+1))\right\}\nonumber\\
\geq& \min\limits_{u_i(k)} \Bigg\{\left(1+ \theta_i\frac{\alpha_i-1}{(1+\theta_i)} \right)r_i(e_i(k), u_i, u_{(j)}) + \left(\alpha_i-\frac{\alpha_i-1}{(1+\theta_i)}\right)J^*_i(e_i(k+1))\Bigg\}\nonumber\\
=&\left(1+ \frac{\alpha_i-1}{(1+\theta_i^{-1})}\right)\min\limits_{u_i(k)} \left\{r_i(e_i(k), u_i, u_{(j)}) + J^*_i(e_i(k+1))\right\}\nonumber\\
=&\left(1+ \frac{\alpha_i-1}{(1+\theta_i^{-1})}\right) J^*_i(e_i(k))
\end{align}

Assume that the inequality (\ref{E36}) holds for $s-1$. Then, we have
\begin{align*}\label{E37}
V_i^s(e_i(k))=&\min\limits_{u_i(k)} \left\{r_i(e_i(k), u_i, u_{(j)}) + V^{s-1}_i(e_i(k+1))\right\}\nonumber\\
\geq& \min\limits_{u_i(k)} \big\{r_i(e_i(k), u_i, u_{(j)}) +\left(1+ \frac{\alpha_i-1}{(1+\theta_i^{-1})^{s-1}}\right)   J^*_i(e_i(k+1))\big\}\nonumber\\
\geq& \min\limits_{u_i(k)} \Bigg\{\left(1+ \frac{(\alpha_i-1)\theta_i^s}{(1+\theta_i)^s} \right)r_i(e_i(k), u_i, u_{(j)}) \nonumber\\
 &+ \left(1+\frac{\alpha_i-1}{(1+\theta_i^{-1})^{s-1}}- \frac{(\alpha_i-1)\theta_i^{s-1}}{(\theta_i+1)^{s}} \right)J^*_i(e_i(k+1))\Bigg\}\nonumber\\
=&\left(1+ \frac{(\alpha_i-1)\theta_i^{s}}{(\theta_i+1)^{s}} \right)\min\limits_{u_i(k)} \left\{r_i(e_i(k), u_i, u_{(j)}) + J^*_i(e_i(k+1))\right\} \nonumber\\
=&\left(1+ \frac{\alpha_i-1}{(1+\theta_i^{-1})^s}\right) J^*_i(e_i(k)).
\end{align*}

The right part of inequality \eqref{20171213Eq-1}
\begin{eqnarray*}\label{E38}
V_i^s(e_i(k))\leq \left(1+ \frac{\beta_i-1}{(1+\theta_i^{-1})^s} \right) J^*_i(e_i(k))
\end{eqnarray*}
can also be proven by the similar development.

Hereinafter, the uniform convergence of value function is demonstrated as the iteration index $s \to \infty$.
When $s \rightarrow \infty$, for $0 < \theta_i < \infty$, we have
\begin{eqnarray*}\label{E39}
\lim \limits_{s\rightarrow \infty} \left(1+ \frac{\alpha_i-1}{(1+\theta_i^{-1})^s} \right) J^*_i(e_i(k))=J^*_i(e_i(k))
\end{eqnarray*}
and
\begin{eqnarray*}\label{E40}
\lim \limits_{s\rightarrow \infty} \left(1+ \frac{\beta_i-1}{(1+\theta_i^{-1})^s} \right) J^*_i(e_i(k))=J^*_i(e_i(k)).
\end{eqnarray*}

Define $V_i^{\infty}(e_i(k)) = \lim \limits_{s\rightarrow \infty} V_i^s(e_i(k))$. Then we can obtain
\begin{eqnarray*}\label{E41}
V_i^\infty(e_i(k)) =J_i^*(e_i(k))
\end{eqnarray*}
Therefore, $V_i(e_i(k))$ converges to $J_i^*(e_i(k))$ as $s \rightarrow \infty$.




\ifCLASSOPTIONcaptionsoff
  \newpage
\fi

\twocolumn



\bibliographystyle{IEEEtran}
\bibliography{mydatabase}

\begin{thebibliography}{10}
\providecommand{\url}[1]{#1}
\csname url@samestyle\endcsname
\providecommand{\newblock}{\relax}
\providecommand{\bibinfo}[2]{#2}
\providecommand{\BIBentrySTDinterwordspacing}{\spaceskip=0pt\relax}
\providecommand{\BIBentryALTinterwordstretchfactor}{4}
\providecommand{\BIBentryALTinterwordspacing}{\spaceskip=\fontdimen2\font plus
\BIBentryALTinterwordstretchfactor\fontdimen3\font minus
  \fontdimen4\font\relax}
\providecommand{\BIBforeignlanguage}[2]{{%
\expandafter\ifx\csname l@#1\endcsname\relax
\typeout{** WARNING: IEEEtran.bst: No hyphenation pattern has been}%
\typeout{** loaded for the language `#1'. Using the pattern for}%
\typeout{** the default language instead.}%
\else
\language=\csname l@#1\endcsname
\fi
#2}}
\providecommand{\BIBdecl}{\relax}
\BIBdecl

\bibitem{cyrc}
Y.~Cao, W.~Yu, W.~Ren, and G.~Chen, ``An overview of recent progress in the
  study of distributed multi-agent coordination,'' \emph{IEEE Transactions on
  Industrial Informatics}, vol.~9, no.~1, pp. 427--438, 2013.

\bibitem{kcm}
S.~Knorn, Z.~Chen, and R.~H. Middleton, ``Overview: Collective control of
  multiagent systems,'' \emph{IEEE Transactions on Control of Network Systems},
  vol.~3, no.~4, pp. 334--347, 2016.

\bibitem{q}
J.~Qin, Q.~Ma, Y.~Shi, and L.~Wang, ``Recent advances in consensus of
  multi-agent systems: A brief survey,'' \emph{IEEE Transactions on Industrial
  Electronics}, vol.~64, no.~4, pp. 4972 -- 4983, 2017.

\bibitem{Oha2015Auto-survey}
K.-K. Oha, M.-C. Park, and H.-S. Ahn, ``A survey of multi-agent formation
  control,'' \emph{Automatica}, vol.~53, pp. 424--440, 2015.

\bibitem{ren2008Book-distributed}
W.~Ren and R.~Beard, \emph{{Distributed Consensus in Multi-Vehicle Cooperative
  Control: Theory and Applications}}.\hskip 1em plus 0.5em minus 0.4em\relax
  London: Springer-Verlag, 2008.

\bibitem{Qu2009Book-cooperative}
Z.~Qu, \emph{{Cooperative Control of Dynamical Systems: Applications to
  Autonomous Vehicles}}.\hskip 1em plus 0.5em minus 0.4em\relax London:
  Springer-Verlag, 2009.

\bibitem{LewisZhang2014Book-cooperative}
F.~Lewis, H.~Zhang, K.~Hengster-Movric, and A.~Das, \emph{{Cooperative Control
  of Multi-Agent Systems: Optimal and Adaptive Design Approaches}}.\hskip 1em
  plus 0.5em minus 0.4em\relax London: Springer-Verlag, 2014.

\bibitem{zld}
H.~Zhang, F.~L. Lewis, and A.~Das, ``Optimal design for synchronization of
  cooperative systems: State feedback, observer and output feedback,''
  \emph{IEEE Transactions on Automatic Control}, vol.~56, no.~8, pp.
  1948--1952, 2011.

\bibitem{hylx}
K.~Hengster-Movric, K.~Y. You, F.~L. Lewis, and L.~Xie, ``Synchronization of
  discrete-time multi-agent systems on graphs using riccati design,''
  \emph{Automatica}, vol.~49, no.~2, pp. 414--423, 2013.

\bibitem{Vamvoudakis2012}
K.~G. Vamvoudakis, F.~L. Lewis, and G.~R. Hudas, ``Multi-agent differential
  graphical games: Online adaptive learning solution for synchronization with
  optimality,'' \emph{Automatica}, vol.~48, no.~8, pp. 1598--1611, 2012.

\bibitem{zlwf}
H.~Zhang, H.~Liang, Z.~Wang, and T.~Feng, ``Optimal output regulation for
  heterogeneous multiagent systems via adaptive dynamic programming,''
  \emph{IEEE Transactions on Neural Networks and Learning Systems}, vol.~28,
  no.~1, pp. 18--29, 2017.

\bibitem{Vamvoudakis2011}
K.~G. Vamvoudakis and F.~L. Lewis, ``Multi-player non-zero-sum games: Online
  adaptive learning solution of coupled {Hamilton-Jacobi} equations,''
  \emph{Automatica}, vol.~47, no.~8, pp. 1556--1569, 2011.

\bibitem{mmns2018}
M.~Mazouchi, M.~B. Naghibi-Sistani, and S.~K.~H. Sani, ``A novel distributed
  optimal adaptive control algorithm for nonlinear multi-agent differential
  graphical games,'' \emph{IEEE/CAA Journal of Automatica Sinica}, vol.~5,
  no.~1, pp. 331--341, 2018.

\bibitem{ztll}
H.~Zhang, T.~Feng, H.~Liang, and Y.~Luo, ``{LQR}-based optimal distributed
  cooperative design for linear discrete-time multiagent systems,'' \emph{IEEE
  Transactions on Neural Networks and Learning Systems}, vol.~28, no.~3, pp.
  599--611, 2017.

\bibitem{jmslv}
Q.~Jiao, H.~Modares, S.~Y. Xu, F.~L. Lewis, and K.~G. Vamvoudakis,
  ``Multi-agent zero-sum differential graphical games for disturbance rejection
  in distributed control,'' \emph{Automatica}, vol.~69, pp. 24--34, 2016.

\bibitem{jj}
Y.~Jiang and Z.~P. Jiang, ``Computational adaptive optimal control for
  continuous-time linear systems with completely unknown dynamics,''
  \emph{Automatica}, vol.~48, no.~10, pp. 2699--2704, 2012.

\bibitem{zl}
H.~Zhang and F.~L. Lewis, ``Adaptive cooperative tracking control of
  higher-order nonlinear systems with unknown dynamics,'' \emph{Automatica},
  vol.~48, no.~7, pp. 1432--1439, 2012.

\bibitem{zzf}
J.~Zhang, H.~Zhang, and T.~Feng, ``Distributed optimal consensus control for
  nonlinear multiagent system with unknown dynamic,'' \emph{IEEE Transactions
  on Neural Networks and Learning Systems}, in press, DOI:
  10.1109/TNNLS.2017.2728622, 2017.

\bibitem{Lewis2011}
F.~L. Lewis and K.~G. Vamvoudakis, ``Reinforcement learning for partially
  observable dynamic processes: Adaptive dynamic programming using measured
  output data,'' \emph{IEEE Transactions on Systems, Man, and Cybernetics, Part
  B: Cybernetics}, vol.~41, no.~1, pp. 14--25, 2011.

\bibitem{Tamimi2007}
A.~Al-Tamimi, F.~L. Lewis, and M.~Abu-Khalaf, ``Model-free {$Q$-Learning}
  designs for discrete-time zero-sum games with application to hinfinity
  control,'' \emph{Automatica}, vol.~47, no.~1, pp. 207--214, 2011.

\bibitem{Liu2015}
D.~Liu, H.~Li, and D.~Wang, ``Error bounds of adaptive dynamic programming
  algorithms for solving undiscounted optimal control problems,'' \emph{IEEE
  Transactions on Neural Networks and Learning Systems}, vol.~26, no.~6, pp.
  1323--1334, 2015.

\bibitem{gjj}
W.~Gao, Y.~Jiang, Z.~P. Jiang, and T.~Chai, ``Output-feedback adaptive optimal
  control of interconnected systems based on robust adaptive dynamic
  programming,'' \emph{Automatica}, vol.~72, pp. 37--45, 2016.

\bibitem{jilie2014}
H.~Zhang, J.~Zhang, G.~H. Yang, and Y.~Luo, ``Leader-based optimal coordination
  control for the consensus problem of multiagent differential games via fuzzy
  adaptive dynamic programming,'' \emph{IEEE Transactions on Fuzzy Systems},
  vol.~23, no.~1, pp. 152--163, 2015.

\bibitem{gu2012}
G.~Gu, \emph{Discrete-Time Linear Systems: Theory and Design with
  Applications}.\hskip 1em plus 0.5em minus 0.4em\relax London:
  Springer-Verlag, 2012.

\bibitem{Wei2011}
Q.~Wei and D.~Liu, ``Nonlinear multi-person zero-sum differential games using
  iterative adaptive dynamic programming,'' in \emph{Proceedings of the 30th
  Chinese Control Conference}, 2011, pp. 2456--2461.

\bibitem{Zhang2013}
H.~Zhang, D.~Liu, and W.~D. Y.~Luo, \emph{Adaptive Dynamic Programming for
  Control:Algorithms and Stability}.\hskip 1em plus 0.5em minus 0.4em\relax
  London: Springer-Verlag, 2013.

\bibitem{Murry2002}
J.~J. Murray, C.~J. Cox, G.~G. Lendaris, and R.~Saeks, ``Adaptive dynamic
  programming,'' \emph{IEEE Transactions on Systems, Man, and Cybernetics, Part
  C: Applications and Reviews}, vol.~32, no.~2, pp. 140--153, 2002.

\end{thebibliography}
%




%





\end{document}